\documentclass[12pt]{article}
\usepackage{amsmath,amssymb}
\usepackage{amsthm}
\usepackage{array}

\newtheorem{prop}{Proposition}

\newtheorem{lemma}{Lemma}

\theoremstyle{note}
\theoremstyle{definition}

\def\rit#1{{\mbox{\rm #1}}\hspace{1mm}}

\def\itemx#1{\item[{\rm(#1)}]}

\def\R{\mathfrak R}
\def\K{\mathfrak K}

\def\C{\mathfrak C}

\def\SS{\mathfrak S}

\def\rit#1{{\mbox{\rm #1}}}

\def\modx#1#2{\equiv #1 \mod ~#2}
\def\nmod#1#2{\not\equiv #1 \mod ~#2}
\def\itemx#1{\item[{\rm(#1)}]}

\def\cusp#1#2{\begin{pmatrix}#1\\ #2\end{pmatrix}}
\def\cuspv#1#2{\begin{bmatrix}#1\\ #2\end{bmatrix}}
\begin{document}
\title{Representation of modular invariant function by generators of a modular function field}
\maketitle
\begin{center}
 N{\sc oburo} I{\sc shii}\footnote[1]{Partially supported by Grand-in-Aid for Scientific Research  No.15540042. 2000 {\it Mathematics Subject Classification}~11G05,14H52}\vspace{5mm}\end{center}

\section{Introduction}
 For a positive integer $N$, let $\Gamma_0(N)$ be a subgroup of $\rit{SL}_2(\mathbf Z)$ defined by
\[
\Gamma_0(N)=\left\{\left. \begin{pmatrix} a & b \\ c & d \end{pmatrix}\in \rit{SL}_2(\mathbf Z)~\right |~ c \equiv 0 \mod N \right\}.
\]
Let $X_0(N)$ be the modular curve associated with $\Gamma_0(N)$ and $g_0(N)$ the genus of $X_0(N)$. It is well-known that a non-cuspidal point $P$ of $X_0(N)$ is corresponding to an isomorphism class of pairs of an elliptic curve $E$ and a cyclic subgroup of $E$ of order $N$. Our problem is to provide a method to determine the elliptic curve $E$ corresponding to $P$, up to isomorphisms, in other words, to compute the $j$-invariant $j(E)$ of $E$.  We define the $j$-invariant $j(P)$ of $P$ by $j(P)=j(E)$. This problem was studied by R.Fricke [1] in the case the genus of $X_0(N)$ is $0$ and by D.Elkies in the case $N$ is a prime number other than $37$ and $X_0(N)$ is elliptic or hyperelliptic. Further the case $X_0(N)$ is hyperelliptic and $N$ is not prime was studied by T.Hibino and N.Murabayashi [2]. Recently the case $X_0(N)$ is bielliptic is studied by T.Yamauchi [6]. In the case $X_0(N)$ is hyperelliptic or bielliptic, their argument is based on a fact that $X_0(N)$ admits a hyperelliptic involution or a bielliptic involution, and they require computation of Fourier coefficients of cusp forms of weight $2$ with respect to $\Gamma_0(N)$ by using the trace formula. Our method used here is rather classical and can be applied to the case $X_0(N)$ is neither hyperelliptic nor bielliptic. 
 Let $A_0(N)$ be the modular function field of $\Gamma_0(N)$. If $g_0(N)=0$, then the modular invariant function $J$ is a rational function of a generator $X_N$ of $A_0(N)$. In the case $g_0(N)\geq 1$, let $X_N$ and $Y_N$ be two generators of $A_0(N)$. Further let $F_N(X,Y)$ be a polynomial such that $F_N(X_N,Y)=0$ is a monic minimal equation of $Y_N$ over $\mathbf C(X_N)$ and let $J=R_N(X_N,Y_N)$ be a representation of $J$ by $X_N$ and $Y_N$. Then $F_N(X,Y)=0$ gives a singular model $\C_N$ of $X_0(N)$. Let $P=(a,b)$ be a non-singular point of $\C_N$. If $R_N(X,Y)$ is defined at $P$, then a value $j=R_N(a,b)$ is the $j$-invariant of $P$. In this article, we shall give a method to construct generators of $A_0(N)$ and to compute a representation of $J$ by the generators.  By our method, we compute $X_N,Y_N,F_N(X,Y)$ and $R_N(X,Y)$ for every $N$ in the range $6\leq N\leq 50$ and also for $N=52$. We remark that $N=52$ is the smallest integer such that $X_0(N)$ is neither hyperelliptic nor bielliptic. Further it is noted that we have in this range all $N$ such that $g_0(N)=0,1$. In section 1, we explain theoretical background of our method. In section 2, we explain how to construct functions needed to compute $F_N(X,Y)$ and $R_N(X,Y)$. In section 3, we give some examples and tables of our computational results. 
\section{Representation of a function by generators}
Let $\R$ be a Riemann surface of genus $g$ and $\K$ the function field of $\R$.  For a point $P$ on $\R$ we denote by $\K(P)$ a subring of $\K$ consisting of all functions regular at any points other than $P$. For $f\in\K(P)$, we shall denote by $d_P(f)$ the order of pole of $f$ at $P$. Assume that $P$ is not a Weierstrass point. Then for every positive integer $n$, there exists a function $F_n\in\K(P)$ such that $d_P(F_n)=n+g$. For example, see \S 6 of Chapter 2 of [5]. We assume that the leading coefficient of the expansion of $F_n$ at $P$ by the local parameter given in advance  is equal to $1$.   
\begin{prop}\label{prop1} Assume that $P$ is not a Weierstrass point. Let \newline
 $F_1,F_2,\dots,F_g,F_{g+1}$ be the functions given as above.  
\begin{enumerate}
\itemx 1 If $g=0$, then $\K$ is generated by $F_1$ over $\mathbf C$.
\itemx 2 If $g>0$, then $\K$ is generated by $F_1$ and $F_2$ over $\mathbf C$.
\itemx 3 If $F\in\K(P)$, then $F$ is a polynomial of $F_1,F_2,\dots,F_g,F_{g+1}$ over $\mathbf C$.
\end{enumerate}
\end{prop} 
\begin{proof} For a function $H\in\K$, we denote by $d(H)$ total degree of poles of $H$. It is well-known that $d(H)$ is equal to the degree $[\K:\mathbf C(H)]$ of $\K$ over a subfield $\mathbf C(H)$. If $g=0$, then $[\K:\mathbf C(F_1)]=d(F_1)=1$. Therefore $\K=\mathbf C(F_1)$. Let $g>0$. Since $d(F_1)=g+1$ and $d(F_2)=g+2$, they are coprime. Therefore $[\K:\mathbf C(F_1,F_2)]=1$. This shows $\K=\mathbf C(F_1,F_2)$. Since $P$ is not a Weierstrass point, for a non-constant function $H\in\K(P)$, we have $d_P(H)>g$ . We shall prove (3) by induction concerning the order $n=d_P(F)$. Divide $n$ by $2g+1$ and let
$n=(2g+1)\ell + k$, where $\ell,k$ are non-negative integers and $k\leq 2g$. Put\[
U=
\begin{cases}F_{g+1}^\ell\qquad &\text{if }k=0,\\
             F_{g+1}^{\ell-1} F_1F_k \qquad &\text{if }k<g+1,\\
             F_{g+1}^\ell F_{k-g} \qquad &\text{if }k\geq g+1

\end{cases}
\]
Then $U\in\K(P)$ and $d_P(U)=n$. Let $c(F)$ be the leading coefficient of the expansion of $F$ by the local parameter at $P$. Then we have $d_P(F-c(F)U)<n$. By the induction hypothesis, $F-c(F)U$ is a polynomial of $F_1,\dots,F_g$ and $F_{g+1}$.
\end{proof}
In the following, we assume $g>1$. Since $d_P(F_1F_{i+2}-F_2F_{i+1})<2g+i+3$, we know, for $i=1,2,\dots,g-1$, $F_1F_{i+2}-F_2F_{i+1}$ is a linear combination over $\mathbf C$ of $F_{i+1}F_1,F_{i}F_1,\dots,F_2F_1,F_1^2,F_{g+1},F_g,\dots,F_2,F_1,1$ as follows:
\begin{equation}\label{eq0}
F_1F_{i+2}-F_2F_{i+1}+\sum_{k=1}^{i+1}a_{i,k}F_{k}F_1+\sum_{k=1}^{g+1}b_{i,k}F_{k}+c_i=0,
\end{equation}
where $a_{i,k},b_{i,k}$ and $c_i$ are complex numbers. From this system of linear relations, we have a representation of $F_i$ by $F_1$ and $F_2$. In fact, let us consider the following system of linear equations of indeterminates $X_1,\dots,X_{g-1}$ with coefficients in $\mathbf C[F_1,F_2]$: 
\newpage
\begin{equation}\label{eq1}
\begin{cases}
&(F_1+b_{1,3})X_{1}+\sum_{k=2}^{g-1}b_{1,k+2}X_k\\
&\phantom{aaaaaaaaaaa}=-a_{1,1}F_1^2-a_{1,2}F_2F_1-b_{1,1}F_1-b_{1,2}F_2-c_1,\\
&(a_{2,3}F_1+b_{2,3}-F_2)X_{1}+(F_1+b_{2,4})X_{2}+\sum_{k=3}^{g-1}b_{2,k+2}X_k\\
&\phantom{aaaaaaaaaaa}=-a_{2,1}F_1^2-a_{2,2}F_2F_1-b_{2,1}F_1-b_{2,2}F_2-c_2,\\
&\sum_{k=1}^{i-2}(a_{i,k+2}F_1+b_{i,k+2})X_k+(a_{i,i+1}F_1+b_{i,i+1}-F_2)X_{i-1}\\
&\phantom{aaaaaaaaaaa}+(F_1+b_{i,i+2})X_{i}+\sum_{k=i+1}^{g-1}b_{i,k+2}X_k\\
&\phantom{aaaaaaaaaaa}=-a_{i,1}F_1^2-a_{i,2}F_2F_1-b_{i,1}F_1-b_{i,2}F_2-c_i,\\
&\phantom{aaaaaaaaaaaaaaaaaaaaaaaaaaaaa}(i=3,\cdots,g-1).
\end{cases}
\end{equation}
Obviously, $\{X_1=F_3,X_2=F_4,\cdots,X_{g-1}=F_{g+1}\}$ is a solution of \eqref{eq1}.
\begin{lemma}\label{lem1} The coefficient matrix $A$ of \eqref{eq1} is regular.
\end{lemma}
\begin{proof} Let $u_{i,j}$ be the $(i,j)$-component of $A$. Then 
\begin{equation}\label{eq2}
u_{i,j}=
\begin{cases}
a_{i,j+2}F_1+b_{i,j+2} \qquad &\text{if }j<i-1,\\
a_{i,i+1}F_1+b_{i,i+1}-F_2 \qquad &\text{if }j=i-1,\\
F_1+b_{i,i+2}\qquad &\text{if }j=i,\\
b_{i,j+2}\qquad &\text{if }j>i.
\end{cases}
\end{equation}
Consider the determinant of $A$:
\[
\left|A\right|=\sum_{\sigma\in S_{g-1}}sgn(\sigma)A(\sigma),
\]
where $S_{g-1}$ is the symmetric group of degree $g-1$ and $sgn(\sigma)$ denotes the signature of $\sigma$ and $A(\sigma)=u_{1,\sigma(1)}\dots u_{g-1,\sigma (g-1)}$. If $\sigma$ is the identity, then 
\[
d_P(A(\sigma))=d_P(\prod_{k=1}^{g-1}(F_1+b_{k,k+2})=g^2-1.
\]
In the case $\sigma$ is not the identity, we shall show $d_P(A(\sigma))<g^2-1$. Suppose there exist $m$ functions $u_{k,\sigma(k)}$ such that $d_P(u_{k,\sigma(k)})=g+1$ and $n$ functions $u_{j,\sigma(j)}$ such that $d_P(u_{j,\sigma(j)})=g+2$ among $u_{1,\sigma(1)},\dots ,u_{g-1,\sigma (g-1)}$. If $m+n=g-1$, then $u_{i,\sigma(i)}$ is not constant for each $i$. By \eqref{eq2}, this shows $\sigma(1)=1$. Since $\sigma(2)\ne 1$, $\sigma(2)=2$. By continuing this discussion, we conclude that $\sigma$ is the identity. This is a contradiction. Therefore $m+n<g-1$. Further we have
\[
\begin{split}
d_P(A(\sigma))&=m(g+1)+n(g+2)=(m+n)(g+1)+n\\ 
&\leq g^2-1+(n-g-1)< g^2-1.
\end{split}
\] 
Hence we have $d_P(|A|)=g^2-1$. In particular, $|A|\ne 0$. 
\end{proof}
By solving \eqref{eq1}, for $i\geq 3$, we obtain a rational function $H_i(X,Y)\in\mathbf C[X,Y]$  such that $F_i=H_i(F_1,F_2)$. For a given function $F\in\K$, the above argument provides us with a method to obtain a rational function $H(X,Y)$ such that $F=H(F_1,F_2)$. Assume $F$ has poles at points $P_1,P_2,\cdots,P_s$ other than $P$. Find  functions $G_i$ which are polynomials of $F_1,F_2,\dots ,F_{g+1}$ so that $F\prod_{k=1}^s(G_i-G_i(P_k))\in \K(P)$. Then Proposition~\ref{prop1} shows $F\prod_{k=1}^s(G_i-G_i(P_k))$ is a polynomial of $F_1,F_2,\dots,F_{g+1}$. By substituting $F_k~(k\geq 3)$ by $H_k(F_1,F_2)$, we obtain a desired rational function $H(X,Y)$. Further, since $F_1$ and $F_2$ have poles only at $P$, by similar arguments in Lemmas 1 and 2 and further in the latter part of Lemma 3 of section 3 of [4], we know the monic minimal equation $F_N(F_1,Y)=0$ of $F_2$ over $\mathbf C(F_1)$ has the following form:
\begin{equation}\label{mineq}
F_N(X,Y)=Y^{g+1}-X^{g+2}+\Phi _{g}(X)Y^g+\cdots +\Phi_1(X)Y+\Phi_0(X),
\end{equation}
where $\Phi_j(X) \in \mathbf C[X]$ and $\deg \Phi_j(X)\leq g+1-j$, for all $j$. 

\section{Modular functions with respect to a group ${\bf \Gamma_0(N)}$}
Let $N$ be a positive integer. We consider two subgroups $\Gamma_0(N)$ and $\Gamma_1(N)$ of $\rit{SL}_2(\mathbf Z)$ defined by
\[
\begin{split}
\Gamma_1(N)&=\left\{\left . \begin{pmatrix} a & b \\ c & d \end{pmatrix}\in \rit{SL}_2(\mathbf Z)~\right |~a \modx 1N,~c \modx 0N \right\},\\
\Gamma_0(N)&=\left\{\left . \begin{pmatrix} a & b \\ c & d \end{pmatrix}\in \rit{SL}_2(\mathbf Z)~\right |~c \modx 0N \right\}.
\end{split}
\]
We denote by $A_1(N)$ and $A_0(N)$ the modular function fields associated with $\Gamma_1(N)$ and $\Gamma_0(N)$ respectively. For a complex number $\tau$ of the complex upper half plane, we denote by $L_{\tau}$ a lattice in $\mathbf C$ generated by $1$ and $\tau$. Let $\wp (z;L_{\tau})$ be the Weierstrass $\wp$-function associated with $L_{\tau}$.  For a vector $\mathbf a=[a_1,a_2,a_3,a_4]$ of integral components such that $a_i\nmod 0N$ for all $i$,~$a_1\nmod {\pm a_2}N$ and $a_3\nmod {\pm a_4}N$, we define a function 
\[
W_{\mathbf a}(\tau)=\frac{\wp (a_1/N;L_{\tau})-\wp (a_2/N;L_{\tau})}{\wp (a_3/N;L_{\tau})-\wp (a_4/N;L_{\tau})}.
\]
 Then the function $W_{\mathbf a}\in A_1(N)$ and it has neither zeros nor poles on the complex upper half plane. For example, see [3]. We note that $W_{\mathbf a}=1$ in the case $(a_1,a_2)\modx {(a_3,a_4)}N$. We shall determine order of $W_{\mathbf a}$ at cusps of $\Gamma_1 (N)$. In [3], all inequivalent cusps of $\Gamma_1 (N)$ are given by pairs of integers $\cusp ut$ such that 
\begin{equation}\label{eq3}
\begin{cases}
\cusp ut ~\text{ for } 1\le t <\frac N2,~1\le u\le D,~(u,D)=1  \text{ and}\vspace{1mm}\\
\cusp ut~\text{ for } t=\frac N2,N,~1\le u\le D/2,~(u,D)=1,
\end{cases}
\end{equation}
where $D=(t,N)$.  Further a local parameter at a cusp $\cusp ut$ is given by $q_D=\exp(2\pi i\tau D/N)$.  For a divisor $D$ of $N$ and an integer $n$, we define two integers $\{n\}_D$ and $\mu_D (n)$  determined uniquely by the following conditions:
\[
\begin{split}
& 0 \le \{n\}_D\le \frac N{2D},\quad \mu_D (n)=\pm 1,\quad n\equiv \mu_D (n)\{n\}_D \mod \frac ND,\\
&\text{if }\{n\}_D=0\text{ or }\frac N{2D}, \text{then }\mu_D (n)=1.
\end{split}
\]
 To simplify the notation, we shall denote $\{n\}_1$ and $\mu_1(n)$ by $\{n\}$ and $\mu (n)$ respectively. We see easily,  for a divisor $D$ of $N$,
\begin{equation}\label{eq4}
\{nD\}=\{n\}_DD\quad \mu (nD)=\mu _D(n).
\end{equation}By similar arguments in Lemmas 1 and 2 of [3], the order of $W_{\mathbf a}$ at a cusp $\cusp ut$ is given as follows.
\begin{prop}\label{prop2}
Let $Q=\cusp ut$ be a cusp of $\Gamma_1(N)$ and $D=(t,N)$. Put $t'=t/D$. Then the order $o_Q(\mathbf a)$ of $W_{\mathbf a}$ at $Q$ is equal to\[\min(\{a_1t'\}_D,~\{a_2t'\}_D)-\min(\{a_3t'\}_D,~\{a_4t'\}_D).\]
\end{prop}
\begin{proof} See Lemmas 1 and 2 of [3] and use ~\eqref{eq4}.
\end{proof}
It is noted that the order is independent of $u$. Next we shall consider cusps of $\Gamma_0(N)$. Let $E_2$ denote the unit matrix of degree $2$. The group $G(N)=\Gamma_0(N)/\{\pm E_2\}\Gamma_1(N)$ operates on a set of inequivalent cusps $\cusp ut$ of the group $\Gamma_1(N)$. We know that every $G(N)$-equivalent class is represented by a cusp $\cusp uD$ given in \eqref{eq3} for a divisor $D$ of $N$ . Further a cusp $\cusp uD$ is $G(N)$-equivalent to a cusp $\cusp vD$  if and only if $u\modx v{(D,N/D)}$. Therefore a cusp $\cusp ut$ is $G(N)$-equivalent to a cusp $\cusp vD$ if and only if $D=(t,N)$ and $tu/D\modx v{(D,N/D)}$. We shall denote by ${\cuspv uD}$ a cusp of $\Gamma_0(N)$ corresponding to the $G(N)$-equivalent class represented by $\cusp uD$. A local parameter at $\cuspv uD$ is given by $q_D^{(D,N/D)}$. Let $\SS\subset \Gamma_0(N)$ be a system of representatives of coset decomposition of $\Gamma_0(N)$ by $\{\pm E_2\}\Gamma_1(N)$. For $F\in A_1(N)$, put $ T(F)=\sum_{A\in\SS}F(A(\tau))$. Then obviously $T(F)\in A_0(N)$. Let $\R(N)$ be a Riemann surface associated with $\Gamma_0(N)$ and let $\K=A_0(N)$. For a cusp $P=\cuspv 11$, we shall construct  functions $F_i\in \K(P)$ with properties described in \S 1 by using functions $T(W_{\mathbf a}W_{\mathbf b})$. These functions have poles only at cusps of $\Gamma_0(N)$. It is noted that a function $T(W_{\mathbf a})$ appears among functions $T(W_{\mathbf a}W_{\mathbf b})$ by taking $\mathbf b=[b_1,b_2,b_1,b_2]$. For $\lambda\in (\mathbf Z/N\mathbf Z)^\times/\{\pm 1\} $, take $M_\lambda\in\Gamma_0(N)$ so that $M_\lambda\modx {\begin{pmatrix}\lambda^{-1} & 0\\0&\lambda\end{pmatrix}}N$. Let $\SS=\{M_\lambda |~\lambda\in (\mathbf Z/N\mathbf Z)^\times/\{\pm 1\}$. By \S 2 of [3], we know $W_{\mathbf a}(M_\lambda(\tau))=W_{\lambda \mathbf a}(\tau)$, where $\lambda \mathbf a$ denotes a scalar multiple of $\mathbf a$ by $\lambda$. Therefore
\[T(W_{\mathbf a}W_{\mathbf b})=\sum_{\lambda}W_{\lambda\mathbf a}W_{\lambda\mathbf b},
\]
where $\lambda$ runs over $(\mathbf Z/N\mathbf Z)^\times/\{\pm 1\}$. Let $P=\cuspv uD$ and denote by $\overline o_P(\mathbf a,\mathbf b)$ the order of $T(W_{\mathbf a}W_{\mathbf b})$ at the cusp $P$. Then by Proposition~\ref{prop2} we see easily
\begin{equation}\label{eqb}
\overline o_P(\mathbf a,\mathbf b)\geq \min_{s}(o_Q(s\mathbf a)+o_Q(s\mathbf b))
\end{equation}
where $Q=\cusp 1D$ and $s$ runs over all integers $s$ such that $(s,N/D)=1,1\leq s \leq N/(2D)$. Since $(u,D)=1$, there exist two integers $c,d$ such that $ud-cD=1$. Let $B(u,D)=\begin{pmatrix} u & c \\ D& d\end{pmatrix}$. Then the $q_D^{(D,N/D)}$-expansion of $T(W_{\mathbf a}W_{\mathbf b})$ at $P$ is defined to be that of $T(W_{\mathbf a}W_{\mathbf b})(B(u,D)(\tau))$, which can be obtained as follows. By definition, we have
\begin{equation}\label{eqa}
\begin{split}
T(W_{\mathbf a}W_{\mathbf b})&(B(u,D)(\tau))=\\
&\sum_{\lambda}\left(\frac{\wp(\frac{\lambda a_1(D\tau+d)}N;L_\tau)-\wp(\frac{\lambda a_2(D\tau+d)}N;L_\tau)}{\wp(\frac{\lambda a_3(D\tau+d)}N;L_\tau)-\wp(\frac{\lambda a_4(D\tau+d)}N;L_\tau)}\right)\times\\
&\phantom{aaaaaaaaa}\left(\frac{\wp(\frac{\lambda b_1(D\tau+d)}N;L_\tau)-\wp(\frac{\lambda b_2(D\tau+d)}N;L_\tau)}{\wp(\frac{\lambda b_3(D\tau+d)}N;L_\tau)-\wp(\frac{\lambda b_4(D\tau+d)}N;L_\tau)}\right).
\end{split}
\end{equation}
Further by Lemma 1 of [3], we know, for integers $r$ and $s$, 
\begin{equation}\label{eq5}
\begin{split}
\wp&\left(\frac{\lambda r(D\tau+d)}N;L_\tau\right)-\wp\left(\frac{\lambda s(D\tau+d)}N;L_\tau\right)=\\
&\sum_{m=0}^\infty\sum_{n=1}^{\infty}n\left[\zeta_d^{n(r\lambda)^*}q_D^{n\{r\lambda\}_D }-\zeta_d^{n(s\lambda)^*}q_D^{n\{s\lambda\}_D }\right]q_D^{mnN/D}\\
&\phantom{aa}+\sum_{m=1}^\infty\sum_{n=1}^\infty n\left[\zeta_d^{-n(r\lambda)^*}q_D^{-n\{r\lambda\}_D }-\zeta_d^{-n(s\lambda)^*}q_D^{-n\{s\lambda\}_D}\right]q_D^{mnN/D},
\end{split}
\end{equation}
where $\zeta_d=\rit{exp}(2\pi id/N)$ and, for an integer $c$, $c^*=\mu_D(c)cd$ .
Therefore by \eqref{eqa} and \eqref{eq5}, we have a desired expansion of $T(W_{\mathbf a}W_{\mathbf b})(B(u,D)(\tau))$. To find a function $F_i$, by \eqref{eqb}, we search for vectors ${\mathbf a}$ and ${\mathbf b}$ so that $\min_{s_1}(o_{Q_1}(s_1\mathbf a)+o_{Q_1}(s_1\mathbf b))=-g-i$ and 
\[\min_{s_D}(o_{Q_D}(s_D\mathbf a)+o_{Q_D}(s_D\mathbf b))\ge -(D,N/D)\]
for any divisor $D\ne 1$ of $N$. Here $Q_t=\cusp 1t$ and $s_t$ runs over all integers $s_t$ such that $(s_t,N/t)=1,1\leq s_t \leq N/(2t)$. If we obtain such vectors ${\mathbf a}$ and ${\mathbf b}$, then we compute $q_1$-expansion of $T(W_{\mathbf a}W_{\mathbf b})((\tau-1)/\tau)$ to see whether $\overline o_{P_1}(\mathbf a,\mathbf b) =-i$. We shall compute $R_N(X,Y)$ in the following procedure.\vspace{3mm}\newline
$\mathbf{1.}$ Let $P=\cuspv 11$ and put $q=q_1$. We search for $F_1,F_2,\dots,F_{g_0(N)+1}$ among functions of the form $\alpha +aT(W_{\mathbf a})+bT(W_{\mathbf b}W_{\mathbf c})+\cdots$.\newline
$\mathbf{2.}$ By (4) and by using $q$-expansions of $F_1$ and $F_2$ at $P$, we compute $F_N(X,Y)$.\newline
$\mathbf{3.}$  For $i=1,2,\dots,g-1$, by (3) of Proposition 1, we express $F_1F_{i+2}-F_2F_{i+1}$ as a linear combination over $\mathbf C$ of $F_{i+1}F_1,F_{i}F_1,\dots,F_2F_1,F_1^2$, $F_{g+1},F_g,$ $\dots,F_2,F_1$ and $1$.\newline
$\mathbf{4.}$ By solving the system of linear equations obtained in $\mathbf 3$,
we have a rational function $H_i(X,Y)$ such that $F_i=H_i(F_1,F_2)$.\newline
$\mathbf {5.}$ For each cusp $Q$  of $X_0(N)$ other than $P$, we search for a function $G_Q$ which is a linear form of $F_1,F_2,\dots,F_{g_0(N)+1}$ and has order of zero at $Q$ as large as possible. \newline
$\mathbf {6.}$ Choose positive integers $m_Q$ so that $J \prod_QG_Q^{m_Q} \in\K(P)$. From (3) of Proposition 1, we obtain 
\[J=\frac{P_N(F_1,\cdots,F_{g_0(N)+1})}{\prod_QG_Q^{m_Q}},\] for a polynomial $P_N(X_1,\dots ,X_{g_0(N)+1})$. If we substitute $F_i$ by $H_i(X,Y)$ in the equation in  $\mathbf 6$  for $i\geq 3$, then we may obtain $R_N(X,Y)$. However, the function $R_N(X,Y)$ will be much complicated as $N$ becomes large. Therefore we would like to present our computational result without substituting $F_i$ by $H_i(F_1,F_2)$.
\section{Examples}
In this section, we give computational results for $N=14$ and $N=52$. The number $N=14$ is the smallest non-prime number such that $X_0(N)$ is an elliptic curve and $N=52$ is the smallest integer such that $X_0(N)$ is neither hyperelliptic nor bielliptic. We use the notation in sections 2 and 3. Furthermore let us denote $q_1$ by $q$ for simplicity. First we treat of the case $N=14$. In this case we have $g_0(14)=1$ and we can take $F_1=T(W_{[5,1,2,1]}),F_2=T(W_{[4,1,3,1]}W_{[5,1,2,1]})$. We know $X_0(14)$ has four cusps $\cuspv 1D,~(D=1,2,7,14)$. Let $F_{i,D}$ be the $q_D$-expansion of $F_i$ at $\cuspv 1D$. Then, by \eqref{eq5}, we have
\begin{equation}\label{eq14}
\begin{cases}
&F_{1,1}=q^{-2}+q^{-1}+2q+2q^2+3q^3+4q^4-2q^5-q^6+q^7\\
&-4q^8-2q^9-6q^{10}-10q^{11}+8q^{12}+6q^{13}-3q^{14}\\
&+10q^{15}+\cdots,\\
&F_{1,2}=-1+8\zeta_{14}^2q_2^2+8\zeta_{14}^3q_2^3+\cdots,\\
&F_{1,7}=7q_7+21q_7^2+\cdots,\\
&F_{1,14}=7+56q_{14}+\cdots,\\
&F_{2,1}=q^{-3}+q^{-2}+3q^{-1}+5+7q+6q^2+5q^3+8q^4+8q^5\\
&-q^6-3q^7-16q^8-11q^9+2q^{10}-26q^{11}-24q^{12}+10q^{13}\\
&+9q^{14}+42q^{15}+\cdots,\\
&F_{2,2}=1+8\zeta_{14} q_2+8\zeta_{14}^3q_2^3+\cdots,\\
&F_{2,7}=-3-21q_7+\cdots,\\
&F_{2,14}=25+224q_{14}+\cdots
\end{cases}
\end{equation}
From two $q$-expansions $F_{1,1}$ and $F_{2,1}$, by \eqref{mineq}, we have an equation between $X=F_1$ and $Y=F_2$:
\[Y^2-X^3+YX-6X^2-Y-18X-12=0\]
Since $J$ has a pole of degree $14/D$ at $\cuspv 1D$, $J(X+1)^4X^2(X-7)$ has a pole only at $\cuspv 11$. By the argument in Proposition~\ref{prop1},we obtain
\[J=\frac{AY+B}{(X+1)^4X^2(X-7)},\]
where
\[
\begin{split}
&A=-7X^{12}-28X^{11}+154X^{10}+1588X^9+5775X^8+11592X^7\\
&+14028X^6+10248X^5+4263X^4+980X^3+4410X^2+196X+49,\\
&B=X^{14}+18X^{13}+62X^{12}-416X^{11}-4665X^{10}-19750X^9\\
&-47712X^8-71184X^7-70977X^6-56762X^5-41850X^4-6672X^3\\
&+5593X^2-882X-196.
\end{split}
\]
The elliptic curve $X_0(14)$ has five $\mathbf Q$-rational points $P=(7,25),[2]P=(0,4),[3]P=(-1,1)$,$[4]P=(0,-3),[5]P=(7,-31)$, except the point $O$ at infinity.   By \eqref{eq14}, we know $O,P,[3]P$ and $[4]P$ correspond to cusps $\cuspv 11$,$\cuspv 1{14}$,\newline $\cuspv 12$ and $\cuspv 17$ respectively. Therefore at these four points, the value of $j$-invariant is $\infty$.  For other two points, by above expression of $J$, we know $j([2]P)=-3375,j([5]P)=16581380$. Next we give results for $N=52$. In this case, $g_0(52)=5$ and we can choose
\[
\begin{split}
&F_1=-\frac{T(W_{[19,3,25,3]}W_{[25,1,19,1]})}2,\quad F_2=\frac {T(W_{[15,2,4,2]}W_{[4,1,15,1]})}3\\
&F_3=-T(W_{[25,1,2,1]}),\quad F_4=T(W_{[23,2,3,2]}W_{[3,1,23,1]}),\\
&F_5=-\frac {T(W_{[16,10,4,10]})}2,\quad F_6=\frac {T(W_{[14,12,13,12]})}2.
\end{split}
\]
We have $q$-expansion $F_{i,1}$ of $F_i$ at the cusp $\cuspv 11$ as follows.
{\scriptsize
\[
\begin{split}
F_{1,1}& = q^{-6}+2q^{-4}+q^{-2}-1+2q^{2}+3q^{6}+2q^{8}+2q^{12}+2q^{14}+2q^{16}+1q^{18}-4q^{22}-4q^{28}\\
&+2q^{30}-6q^{32}-3q^{34}+4q^{36}-2q^{38}-4q^{40}-2q^{42}-q^{46}+8q^{48}+q^{50}+7q^{54}+9q^{58}+4q^{60}\\
&-6q^{62}+2q^{64}+10q^{66}+4q^{68}+q^{70}+2q^{72}-16q^{74}-4q^{76}-q^{78}-12q^{80}-2q^{82}-16q^{84}\\
&-5q^{86}+8q^{88}-4q^{90}-18q^{92}-10q^{94}-5q^{98}+28q^{100}+5q^{102}+2q^{104}+24q^{106}+2q^{108}\\
&+31q^{110}+8q^{112}-16q^{114}+8q^{116}+30q^{118}+16q^{120}-q^{122}+12q^{124}-44q^{126}-6q^{128}\\
&-3q^{130}-40q^{132}-4q^{134}-52q^{136}-10q^{138}+30q^{140}-12q^{142}-54q^{144}-22q^{146}-22q^{150}\\
&+72q^{152}+10q^{154}+6q^{156}+60q^{158}+12q^{160}+80q^{162}+8q^{164}-48q^{166}+16q^{168}\\
&+90q^{170}+34q^{172}-3q^{174}+40q^{176}-235/2q^{178}-43/2q^{180}+\cdots,
\end{split}
\]

\[
\begin{split}
F_{2,1}&= q^{-7}+1/3q^{-5}+2/3q^{-3}+4/3q^{-1}+2+2/3q+2/3q^{3}+1/3q^{5}+5/3q^{7}+2q^{9}+7/3q^{11}\\
&+1/3q^{13}+5/3q^{15}+2q^{17}-7/3q^{19}+q^{21}-2/3q^{23}-4q^{25}+4/3q^{29}+4/3q^{31}-10/3q^{33}\\
&-4q^{35}-5q^{37}-q^{39}-q^{41}-14/3q^{43}+20/3q^{45}-q^{47}+4/3q^{49}+34/3q^{51}+2/3q^{53}\\
&-14/3q^{55}-5/3q^{57}+8q^{59}+22/3q^{61}+28/3q^{63}+2q^{65}-1/3q^{67}+34/3q^{69}-13q^{71}+q^{73}\\
&-2/3q^{75}-22q^{77}-2q^{79}+26/3q^{81}+11/3q^{83}-55/3q^{85}-14q^{87}-59/3q^{89}-4q^{91}-2q^{93}\\
&-68/3q^{95}+73/3q^{97}-10/3q^{99}+10/3q^{101}+122/3q^{103}+10/3q^{105}-16q^{107}-1q^{109}\\
&+37q^{111}+86/3q^{113}+118/3q^{115}+8q^{117}+17/3q^{119}+128/3q^{121}-46q^{123}+7q^{125}\\
&-26/3q^{127}-230/3q^{129}-4q^{131}+34q^{133}+67/3q^{135}-200/3q^{137}-160/3q^{139}-199/3q^{141}\\
&-44/3q^{143}-5q^{145}-78q^{147}+244/3q^{149}-31/3q^{151}+34/3q^{153}+406/3q^{155}+6q^{157}\\
&-190/3q^{159}-37q^{161}+344/3q^{163}+260/3q^{165}+332/3q^{167}+25q^{169}+19/3q^{171}+134q^{173}\\
&-138q^{175}+44/3q^{177}-58/3q^{179}-710/3q^{181}+\cdots
\end{split}
\]
\[
\begin{split}
F_{3,1}&= q^{-8}+q^{-6}+2q^{-4}-15+3q^{2}+2q^{4}+2q^{6}+3q^{8}+2q^{10}+2q^{12}+4q^{14}+2q^{16}-4q^{22}\\
&+5q^{24}+q^{26}-8q^{28}-4q^{30}-2q^{32}-7q^{34}-2q^{36}-10q^{40}+2q^{42}+8q^{48}-10q^{50}\\
&-2q^{52}+18q^{54}+14q^{56}+3q^{58}+16q^{60}+8q^{62}-2q^{64}+22q^{66}-8q^{68}-2q^{70}\\
&+q^{72}-14q^{74}+20q^{76}+4q^{78}-36q^{80}-36q^{82}-4q^{84}-30q^{86}-24q^{88}+2q^{90}-42q^{92}\\
&+12q^{94}+8q^{96}-2q^{98}+26q^{100}-46q^{102}-9q^{104}+74q^{106}+74q^{108}+4q^{110}+56q^{112}\\
&+48q^{114}-4q^{116}+80q^{118}-26q^{120}-11q^{122}+4q^{124}-40q^{126}+94q^{128}+17q^{130}\\
&-136q^{132}-140q^{134}-106q^{138}-90q^{140}+12q^{142}-150q^{144}+64q^{146}+16q^{148}\\
&-8q^{150}+62q^{152}-174q^{154}-30q^{156}+236q^{158}+260q^{160}-11q^{162}+188q^{164}+172q^{166}\\
&-28q^{168}+258q^{170}-126q^{172}-36q^{174}+16q^{176}-100q^{178}+6q^{179}+320q^{180}+18q^{181}+\cdots
\end{split}
\]
\[
\begin{split}
F_{4,1}& =q^{-9}-q^{-8}+q^{-7}+q^{-5}+3q^{-3}+q^{-2}+3q^{-1}+6+3q-q^{2}+3q^{3}-2q^{4}+5q^{5}+q^{6}\\
&+3q^{7}-q^{8}+2q^{9}-2q^{10}+5q^{11}-2q^{14}+5q^{15}+5q^{17}+q^{18}+q^{19}+3q^{21}-4q^{23}\\
&-5q^{24}-6q^{25}-q^{26}-3q^{27}+4q^{28}-q^{29}+6q^{30}-10q^{31}-4q^{32}-4q^{33}+4q^{34}+3q^{35}\\
&+6q^{36}-1q^{37}-2q^{38}+q^{39}+6q^{40}-5q^{41}-4q^{42}-15q^{43}-2q^{45}-q^{46}-3q^{47}\\
&+10q^{49}+1q^{50}+17q^{51}+2q^{52}+8q^{53}-11q^{54}-3q^{55}-14q^{56}+25q^{57}+6q^{58}+16q^{59}\\
&-12q^{60}-9q^{61}-14q^{62}+28q^{63}+4q^{64}-3q^{65}-12q^{66}+5q^{67}+12q^{68}+36q^{69}+3q^{70}\\
&+15q^{71}+q^{72}+q^{73}-2q^{74}-26q^{75}-24q^{76}-32q^{77}-5q^{78}-22q^{79}+24q^{80}+9q^{81}\\
&+34q^{82}-59q^{83}-12q^{84}-43q^{85}+25q^{86}+24q^{87}+32q^{88}-61q^{89}-6q^{90}+7q^{91}+24q^{92}\\
&-12q^{93}-22q^{94}-82q^{95}-8q^{96}-39q^{97}-3q^{98}-2q^{99}+2q^{100}+63q^{101}+51q^{102}+57q^{103}\\
&+11q^{104}+47q^{105}-50q^{106}-16q^{107}-72q^{108}+119q^{109}+27q^{110}+89q^{111}-48q^{112}\\
&-58q^{113}-64q^{114}+120q^{115}+12q^{116}-15q^{117}-50q^{118}+27q^{119}+42q^{120}+173q^{121}\\
&+10q^{122}+82q^{123}+8q^{124}+7q^{125}-4q^{126}-127q^{127}-100q^{128}-105q^{129}-20q^{130}\\
&-83q^{131}+96q^{132}+35q^{133}+136q^{134}-229q^{135}-52q^{136}-170q^{137}+96q^{138}+123q^{139}\\
&+120q^{140}-221q^{141}-24q^{142}+32q^{143}+96q^{144}-47q^{145}-86q^{146}-334q^{147}-16q^{148}\\
&-174q^{149}-14q^{150}-7q^{151}+10q^{152}+238q^{153}+184q^{154}+186q^{155}+36q^{156}+149q^{157}\\
&-176q^{158}-78q^{159}-248q^{160}+437q^{161}+91q^{162}+322q^{163}-180q^{164}-251q^{165}-220q^{166}\\
&+404q^{167}+44q^{168}-64q^{169}-168q^{170}+79q^{171}+160q^{172}+608q^{173}+33q^{174}+340q^{175}\\
&+24q^{176}+7q^{177}-14q^{178}-449q^{179}-340q^{180}-318q^{181}+\cdots\\
F_{5,1}& = q^{-10}-q^{-8}-q^{-4}+3q^{-2}-5-q^{2}-q^{4}+q^{6}+q^{8}+3q^{12}+q^{14}-2q^{16}+q^{18}\\
&+2q^{20}+7q^{22}-13q^{24}-3q^{26}+6q^{28}+5q^{30}+2q^{32}-4q^{34}+q^{36}-13q^{38}-6q^{40}-2q^{42}\\
&+8q^{44}-5q^{46}-16q^{48}+35q^{50}+7q^{52}-13q^{54}-14q^{56}-14q^{58}+14q^{60}-2q^{62}+34q^{64}\\
&+24q^{66}+12q^{68}-23q^{70}+1q^{72}+36q^{74}-86q^{76}-16q^{78}+28q^{80}+38q^{82}+32q^{84}-39q^{86}\\
&-78q^{90}-63q^{92}-23q^{94}+48q^{96}-27q^{98}-93q^{100}+199q^{102}+37q^{104}-62q^{106}-81q^{108}\\
&-69q^{110}+96q^{112}+12q^{114}+174q^{116}+142q^{118}+50q^{120}-110q^{122}+54q^{124}+200q^{126}\\
&-414q^{128}-76q^{130}+116q^{132}+155q^{134}+160q^{136}-212q^{138}-31q^{140}-361q^{142}\\
&-306q^{144}-110q^{146}+242q^{148}-100q^{150}-390q^{152}+812q^{154}+145q^{156}-210q^{158}\\
&-300q^{160}-333q^{162}+434q^{164}+71q^{166}+708q^{168}+626q^{170}+211q^{172}-479q^{174}\\
&+198q^{176}+750q^{178}-1547q^{180}+\cdots 
\end{split}
\]
\[
\begin{split}
F_{6,1}&= q^{-11}+2q^{-9}+4q^{-7}+4q^{-5}+6q^{-3}+8q^{-1}+6+8q+8q^3+13q^5+12q^7+12q^9\\
&+16q^{11}+q^{13}+12q^{15}+16q^{17}+2q^{21}+2q^{23}-12q^{25}-8q^{27}-28q^{31}-24q^{33}-16q^{35}\\
&-33q^{37}-2q^{39}-14q^{41}-44q^{43}+6q^{45}+16q^{47}-12q^{49}+32q^{51}+32q^{53}-12q^{55}+82q^{57}\\
&+76q^{59}+30q^{61}+76q^{63}+6q^{65}+24q^{67}+106q^{69}+8q^{71}-64q^{73}+36q^{75}-42q^{77}-96q^{79}\\
&+32q^{81}-200q^{83}-199q^{85}-58q^{87}-162q^{89}-16q^{91}-67q^{93}-234q^{95}-48q^{97}+156q^{99}\\
&-86q^{101}+28q^{103}+212q^{105}-64q^{107}+425q^{109}+440q^{111}+128q^{113}+312q^{115}+37q^{117}\\
&+164q^{119}+488q^{121}+120q^{123}-346q^{125}+222q^{127}+16q^{129}-416q^{131}+138q^{133}-872q^{135}\\
&-882q^{137}-260q^{139}-551q^{141}-78q^{143}-342q^{145}-948q^{147}-269q^{149}+756q^{151}-528q^{153}\\
&-152q^{155}+836q^{157}-292q^{159}+1746q^{161}+1708q^{163}+460q^{165}+968q^{167}+158q^{169}\\
&+680q^{171}+1742q^{173}+564q^{175}-1560q^{177}+1100q^{179}-10q^{180}+524q^{181}+\cdots
\end{split}
\]
}
By \eqref{mineq} in section 2, we have a relation between $F_1$ and $F_2$:  
\[
\begin{split}
&729F_2^6-8748F_2^5+(10368F_1+8748F_1^2-3888)F_2^4\\
&+(-82944F_1+264384-69984F_1^2)F_2^3\\
&+(1038528F_1^2+15552F_1^4+1645056+2430144F_1+177264F_1^3)F_2^2\\
&-(709056F_1^3+3594240F_1^2+9057024F_1+62208F_1^4+9068544)F_2\\
&+6739136F_1-14067F_1^6-406364F_1^4-729F_1^7+8537280-111148F_1^5\\
&-516848F_1^3+1370992F_1^2=0.
\end{split}
\]

This is a defining equation of $A_0(52)$. Further we obtain a system of linea relations of $F_3,F_4,F_5,F_6$ over $\mathbf C[F_2,F_3]$ from \eqref{eq0}.
\[
\begin{cases}
&(9F_1+170)F3+52F_5+36F_2-9F_2^2-21F_1^2+4F_1+2804=0,\\
&(9F_1-38)F_4+9F_3F_1-15F_1F_2-9F_1^2-9F_2F_3+20-20F_3\\
&\phantom{aaaa}-240F_2+26F_6+140F_1=0,\\
&(9F_1-12)F_5-9F_4F_1-6F_3F_1+24F_1F_2-9F_2F_4-26F_6+46F_1\\
&\phantom{aaaa}+56F_4+27F_1^2+168F_2-1532-70F_3=0,\\
&(9F_1-64)F_6-42F_4F_1-42F_3F_1+48F_1F_2-9F_2F_5+282F_2\\
&\phantom{aaaaa}+18F_5+82F_4+42F_1^2-568F_1+476+82F_3=0.
\end{cases}
\]
By solving them, for $i=3,4,5$ and $6$, we know that $F_i=U_i/\Delta$, where
{\small
\[
\begin{split}
&\Delta=(9828+6318F_1)F_2^2+(-25272F_1-39312)F_2+136608F_1+729F_1^4\\
&\phantom{aaaaa}+12960F_1^3+158512+59360F_1^2,\\
U_3=&2106F_2^4-16848F_2^3+(729F_1^3-200160-1512F_1^2-127620F_1)F_2^2\\
&+(510480F_1+935424+6048F_1^2-2916F_1^3)F_2-2703520\\
&+1701F_1^5+9486F_1^4-740880F_1^2-2013232F_1-112392F_1^3,\\
U_4=&-\{2106F_2^4+(5454F_1+1512-729F_1^2)F_2^3\\
&+(-231696-119628F_1-3456F_1^2+729F_1^3)F_2^2\\
&+(-43524F_1^3-219792F_1+45792-195576F_1^2-2916F_1^4)F_2\\
&+143200F_1+8190F_1^4+13144F_1^3+972F_1^5+8608+39520F_1^2\},\\
U_5=&-\{(-729F_1+5184)F_2^4+(-41472+5832F_1)F_2^3\\
&+(-179064-2673F_1^3-40788F_1^2-151560F_1)F_2^2\\
&+(1048032+10692F_1^3+559584F_1+163152F_1^2)F_2+152576F_1^3\\
&+328880F_1+376816F_1^2-290976+27894F_1^4+2025F_1^5\},\\
U_6=&3\{243F_2^5-2430F_2^4+(2088+6066F_1+2025F_1^2)F_2^3\\
&+(46008-12150F_1^2-23760F_1)F_2^2\\
&+(430880F_1+403408+176824F_1^2+32600F_1^3+2565F_1^4)F_2\\
&-738048-39280F_1^3-592576F_1-202528F_1^2-3672F_1^4\}.
\end{split}
\]
}
Further we know $J(F_1+2)^{13}(F_1+3)(F_1-10)$ has poles only at the cusp $\cuspv 11$ and is a polynomial of $F_6$ of degree $12$ with coefficients in $\mathbf Q[F_1,F_2,\cdots,F_5]$. 
Thus
\[J(F_1+2)^{13}(F_1+3)(F_1-10)=\sum_{i=0}^{12}C_iF_6^i,\] 
where
{\small
\[
\begin{split}
&C_{12}=-28F_1+7F_3+6461/3+F_5,\\
&C_{11}=-8021/3F_4+433F_1F_2-9533/3F_3+125F_1^2+8471/3F_1\\
&-125F_1F_4-72F_5-125F_1F_3-26905F_2-367930/3,\\
\end{split}
\]
\[
\begin{split}
&C_{10}=144932194/3+1775730F_2+176462F_4+20191915/3F_1\\
&+161/3F_1F_5-312308F_5-865108/3F_1^2+197872F_1F_3\\
&-28578F_1F_2+3211861/3F_3+8250F_1F_4,\\
&C_9=18833520F_5-3220F_1F_5-2252927810-525569928F_1\\
&+32134979F_1^2-264445499F_2-5908104F_1F_2+5065804F_4\\
&-26705139F_1F_3-15327819F_1F_4-47918416F_3,\\
&C_8=-241587607/3F_1F_5+68498821957/3F_1+12362268546F_2\\
&+349901856F_1F_2+742849515037/3+1303585160F_5\\
&-464132376F_4+1871964406F_1F_3+40941024352/3F_3\\
&-9194076949/3F_1^2+818792226F_1F_4,\\
&C_7=142425253754F_1^2+24523041006F_1F_2-21515774746F_4\\
&-10557545749788-78860671872F_5-49143435658F_1F_4\\
&-859116229730F_1-653187112730F_3+3868183792F_1F_5\\
&-89865695818F_1F_3-231020246178F_2,\\
&C_6=940986280010659/3+129839421335F_1F_5+1241173308228F_4\\
&+2723618645013F_1F_3+1653156637524F_2+3836730854972F_5\\
&-1269249923340F_1F_2-3892522628187F_1^2+22177550872097F_1\\
&+1514993425764F_1F_4+20226951722173F_3,\\
&C_5=-46452442499126F_1F_3-48537181735174/3F_4\\
&-18106253311940762/3-96733933399152F_5\\
&-1157707902426130/3F_3+32629773607505F_1F_2\\
&+39496097055685F_2+62515859975054F_1^2\\
&-1077076459078412/3F_1-6624064232892F_1F_5\\
&-22475089275218F_1F_4,\\
&C_4=428791323564693F_1F_3+11731631172861482/3F_1\\
&+11184323737178309/3F_3-504606917168286F_1F_2\\
&+168348316204284F_1F_4+1137997144652583F_5\\
&+182185716795601865/3+13458281966236F_4\\
&-1208394806691126F_2-1510986531671354/3F_1^2\\
&+340961275926865/3F_1F_5,\\
\end{split}
\]
\[
\begin{split}
&C_3=11207728140922953F_2-1044316459487816F_1F_5\\
&+4054268582329185F_1F_2-6116228281311912F_5\\
&+875320846743161F_1^2-1694546558285089F_1F_3\\
&-281271228548233F_1F_4-259318888243793602\\
&+931316020753191F_4-14452782305632417F_3\\
&-25842666907378881F_1,\\
&C_2=-12025887423129076F_3+5000846468065935F_1F_5\\
&-34054053314391234F_2+9670685936191636F_5\\
&-5096436215179806F_1F_4+65110227685004338F_1\\
&-136194570485867444/3-4795455072849534F_4\\
&-14187987869301234F_1F_2+13126899329056666F_1^2\\
&-2654467381722550F_1F_3,\\
&C_1=31686984302424144F_5+935824045548092/3F_4\\
&-9637609703523444F_1F_5-25008719190731908F_2\\
&+673556866266135956/3F_3+9146416286935004F_1F_2\\
&+41020836187089844F_1F_3+10355838379726442488/3\\
&+399640736332326652/3F_1-77647177824272164F_1^2\\
&+32635980113943508F_1F_4,\\
&C_0=-82099660820242200F_1F_3+32235149516099208F_1F_2\\
&-95452660752649056F_5-398763269552541040F_3\\
&-6774786532164971264+1525442326328940F_1F_5\\
&-53150472555074088F_1F_4+183431709096990792F_2\\
&+117082538287640184F_1^2+26800586276923736F_4\\
&-635695300913479112F_1.
\end{split}
\]
}
For every $N$ in the range $6\leq N\leq 50$, we have a table of functions $F_1,F_2,\dots,F_{g_0(N)+1}$ and a defining equation $F_N(X,Y)$ and a representation of $J$ by $F_1,F_2$, in almost all cases, via a representation $H_i(X,Y)$ of $F_i$ by $F_1,F_2$.  But it takes too much space to include the table. The essential but difficult process in our computation is that of searching for $F_1,F_2,\dots,F_{g_0(N)+1}$. If these functions are once obtained, then almost automatically $F_N(X,Y)$, $H_i(X,Y)$ and a representation of $J$ by $F_1,\dots,F_{g_0(N)+1}$ can be computed. Therefore, in the case $g_0(N)\geq 3$ or $N\geq 37$, we include only $F_1,\dots,F_{g_0(N)+1}$ in the tables below. The complete table is going to be published on our web site: http://www.las.osakafu-u.ac.jp/$^\sim$ishii/ . In the following tables, for vectors $\mathbf a,\mathbf b,\mathbf c$ and numbers $\alpha,a,b$, the notation $\alpha+a\mathbf a+b\mathbf b*\mathbf c+\cdots$ denotes a function $\alpha +aT(W_{\mathbf a})+bT(W_{\mathbf b}W_{\mathbf c})+\cdots.\vspace{3mm}$\newline
(i) the case $g_0(N)=0$. 
{\scriptsize
\[
\setlength{\extrarowheight}{12pt}
\begin{array}{|l|l|l|}\hline
N&X=F_1&R_N(X)~(J=R_N(F_1))\\ \hline
6&[1,3,2,3]&\frac{(X-3)^3(X^3-9X^2+3X-3)^3}{(X-1)^3X^2(X-9)}\\ \hline
7&[2,1,4,1]&\frac{(X^2-3X+9)(X^2-11X+25)^3}{X-8}\\ \hline
8&[2,1,4,1]&\frac{(X^4-16X^2+16)^3}{(X-4)(X+4)X^2}\\ \hline
9&-3+[2,1,4,1]&\frac{X^3(X^3-24)^3}{X^3-27}\\ \hline
10&-[1,2,4,2]& \frac{(X^6-4X^5+16X+16)^3}{(X-1)^2(X+4)X^5}\\ \hline
12&-3+[5,1,2,1]&\frac{(X^2-3)^3(X^6-9X^4+3X^2-3)^3}{X^4(X-3)(X+3)(X-1)^3(X+1)^3}\\ \hline
13&-5+[5,1,3,1]&\frac{(X^2+X+7)(X^4-X^3+2X^2-9X+3)^3}{X-2}\\ \hline
16& 2+\frac 12[6,1,3,1]& \frac{(X^8-16X^4+16)^3}{X^4(X-2)(X+2)(X^2+4)}\\ 
  &*[7,1,2,1]& \\ \hline
18&-1+[2,1,7,1]&\frac{(X^3-2)^3(X^9-6X^6-12X^3-8)^3}{X^9(X-2)(X^2+2X+4)(X+1)^2(X^2-X+1)^2}\\ \hline
25&-4+[10,2,7,2]&\frac{(X^{10}+10X^8+35X^6-12X^5+50X^4-60X^3+25X^2-60X+16)^3}{(X-1)(X^4+X^3+6X^2+6X+11)}\\ \hline
\end{array}
\]
}
\noindent
(ii) the case $g_0(N)=1$. In this case, we can write $R_N(X,Y)=(A+BY)/C$ with polynomials $A,B,C$ of $X$. The table has the following format.
\[
\begin{array}{|l|c|}\hline
N&X=F_1,~Y=F_2\\ \cline{2-2}
 &F_N(X,Y)=0\\ \cline{2-2}
 &A \\\cline{2-2}
  &B\\\cline{2-2}
   &C\\\hline
\end{array}
\]
{\scriptsize
\[
\renewcommand{\arraystretch}{1.4}
\begin{array}{|l|l|}\hline
11&[2,1,5,1],\quad [2,1,3,1]\\ \cline{2-2}
 &Y^2-5Y-X^3+7X^2-6X+18=0\\ \cline{2-2}
 &-(11X^6-278X^5+1523X^4-1514X^3+974X^2-11777X+12992) \\\cline{2-2}
  &X^5+13X^4-841X^3+5685X^2-10974X+6049\\\cline{2-2}
   &X-18\\\hline
14&[5,1,2,1],\quad [4,1,3,1]*[5,1,2,1] \\\cline{2-2}
  &Y^2-X^3+YX-6X^2-Y-18X-12=0\\\cline{2-2}
  &X^{14}+18X^{13}+62X^{12}-416X^{11}-4665X^{10}-19750X^9-47712X^8-71184X^7\\
  &-70977X^6-56762X^5-41850X^4-6672X^3+5593X^2-882X-196\\\cline{2-2}
  &-(7X^{12}+28X^{11}-154X^{10}-1588X^9-5775X^8-11592X^7-14028X^6\\
  &-10248X^5-4263X^4-980X^3-4410X^2-196X-49)\\\cline{2-2}
  &(X+1)^4X^2(X-7)\\\hline
15&[7,1,3,1],\quad [2,1,7,1] \\\cline{2-2}
  &Y^2-X^3-YX+5X^2-Y+3X-6=0  \\\cline{2-2}
  &-8X^{14}+172X^{13}-1309X^{12}+4682X^{11}-8857X^{10}+10150X^{9}-8004X^{8}\\
  &+1860X^{7}+1110X^{6}+8150X^{5}-11425X^{4}+1375X^{3}+2875X^{2}-2000X+500\\\cline{2-2}
  &X^{13}-5X^{12}-180X^{11}+1925X^{10}-7525X^{9}+13878X^{8}-11370X^{7}\\
  &+330X^{6}+6075X^{5}-3775X^{4}+875X^{3}+1125X^{2}-875X+250 \\\cline{2-2}
  &(X-10)X^3(X-1)^3\\\hline
17&-9+[5,1,2,1],\quad [2,1,4,1]  \\\cline{2-2} 
  & Y^2+XY-12Y-X^3-2X^2-7X+50=0  \\\cline{2-2}
  &-8X^9-51X^8+2X^7+1561X^6+9036X^5+18089X^4+3890X^3-39793X^2\\
  &-37016X-766\\\cline{2-2}
  &X^8+29X^7+12X^6-797X^5-2848X^4-2877X^3+3176X^2+6109X+1291 \\\cline{2-2}
  &X-6\\\hline
19&8+[7,1,1,2],\quad [2,1,4,1] \\\cline{2-2} 
  &Y^2-X^3-3X^2-7X+52+2YX-11Y=0\\\cline{2-2}
  &X^{10}-21X^9+49X^8+341X^7-1202X^6-1145X^5+8429X^4-2519X^3-17952X^2\\
  &+10677X+2036\\\cline{2-2}
  &X^9-6X^8-3X^7+106X^6-182X^5-607X^4+1281X^3+850X^2-1972X+113\\\cline{2-2}
  &X-4\\\hline
20&\frac 12[6,1,8,1],\quad [9,1,2,1] \\\cline{2-2} 
&Y^2-2XY-10Y-X^3+3X^2+5X+25=0\\\cline{2-2}
&(X^6-10X^5+35X^4-60X^3+55X^2-10X+5)^3\\\cline{2-2}
&0\\\cline{2-2}
&(X-1)^5X^2(X-5)\\\hline
\end{array}
\]
\[
\renewcommand{\arraystretch}{1.2}
\begin{array}{|l|l|}\hline
21&-\frac{1}2[5,1,2,1],\quad [6,1,9,1]\\\cline{2-2}
&Y^2-X^3+5XY+6X^2-10Y-21X+26=0\\\cline{2-2}
&-8X^{20}-108X^{19}-290X^{18}+1377X^{17}+7889X^{16}+12618X^{15}+61888X^{14}\\
&+464481X^{13}+1422577X^{12}+1263600X^{11}-3542450X^{10}-12114423X^9\\
&-14689335X^8-3572902X^7+12903306X^6+19284316X^5+13160184X^4\\
&+4506178X^3+197870X^2-400970X-38098\\\cline{2-2}
&X^{19}+61X^{18}+804X^{17}+3677X^{16}-2665X^{15}-96660X^{14}-435876X^{13}\\
&-917535X^{12}-594000X^{11}+1758543X^{10}+5284707X^9+6372696X^8\\
&+2782547X^7-2593906X^6-5046093X^5-3771279X^4-1501035X^3\\
&-252783X^2+19727X+929\\\cline{2-2}
&(X+1)^7(X+2)^2(X-5)\\\hline
24&3-[9,1,2,1],\quad -[4,1,2,1]\\\cline{2-2}
&Y^2+4Y-X^3+3X+4-2X^2=0\\\cline{2-2}
&(X^2-3)^3(X^6-9X^4+3X^2-3)^3\\\cline{2-2}
&0\\\cline{2-2}
&X^4(X-3)(X+3)(X-1)^3(X+1)^3\\\hline
32&2-\frac 14*[12,2,6,2]*[14,2,4,2],\quad [8,2,13,2]\\\cline{2-2}
&Y^2+2XY-14Y-21X+54-X^3+4X^2=0\\\cline{2-2}
&(X^8-16X^4+16)^3\\\cline{2-2}
&0\\\cline{2-2}
&X^4(X-2)(X+2)(X^2+4)\\\hline
36&-1+\frac 12*[4,2,14,2],\quad -\frac 13*[6,9,7,9]\\\cline{2-2}
&Y^2+2Y-X^3=0\\\cline{2-2}
&(X^3-2)^3(X^9-6X^6-12X^3-8)^3\\\cline{2-2}
&0\\\cline{2-2}
&X^9(X-2)(X^2+2X+4)(X+1)^2(X^2-X+1)^2\\\hline
49&-6+[14,4,2,4]*[7,18,3,18],\quad [6,7,21,7]*[21,1,6,1]\\\cline{2-2}
&Y^2-25Y-3XY-X^3+38X+159=0\\\cline{2-2}
&23X^{31}+1989X^{30}+47515X^{29}+425325X^{28}+517993X^{27}-18575564X^{26}\\
&-161213955X^{25}-674629460X^{24}-1844931540X^{23}-5058475096X^{22}\\
&-19271899681X^{21}-67144337491X^{20}-149815471183X^{19}-145903004163X^{18}\\
&+184725017930X^{17}+805084089676X^{16}+950237564249X^{15}-207250599904X^{14}\\
&-1869860337996X^{13}-1778552954157X^{12}+526953944216X^{11}+2216883047056X^{10}\\
&+1145927733576X^9-817610262080X^8-1063185458944X^7-103086959808X^6\\
&+327504300032X^5+122286424064X^4-30972084736X^3-19391012864X^2\\
&-239443968X+438534144\\\cline{2-2}
&X^{30}+299X^{29}+14652X^{28}+281120X^{27}+2799090X^{26}+16541749X^{25}+63288344X^{24}\\
&+179306368X^{23}+511650280X^{22}+1771232319X^{21}+5509411754X^{20}+11262306958X^{19}\\
&+9369410141X^{18}-17331130292X^{17}-63049874238X^{16}-68252518528X^{15}\\
&+25007132730X^{14}+147786037196X^{13}+128818604824X^{12}-53364870160X^{11}\\
&-173818848728X^{10}-79535252304X^9+71127434848X^8+81329280896X^7\\
&+3678398976X^6-26727608320X^5-8668975616X^4+2886603776X^3+1485400064X^2\\
&-20078592X-35979264\\\cline{2-2}
&X^7+14X^4+56X^3+70X^2-28X-113\\\hline
\end{array}
\]  
} 
(iii) the case $g_0(N)=2,N\leq 31$. In this case, we can write $R_N(X,Y)=(A+BY+CY^2)/D$, with polynomials $A,B,C,D$ of $X$. The table has the following format.
\[
\begin{array}{|l|c|}\hline
N&X=F_1,~Y=F_2\\ \cline{2-2}
 &F_N(X,Y)=0\\ \cline{2-2}
 &A \\\cline{2-2}
  &B\\\cline{2-2}
   &C\\\cline{2-2}
   &D\\\hline
\end{array}
\] 
{\scriptsize
\[
\renewcommand{\arraystretch}{1.2}
\begin{array}{|l|l|l|}\hline
22&-3+[8,2,3,2]*[3,1,8,1],\quad 6-[4,2,10,2]*[10,1,4,1]-[8,2,3,2]*[3,1,8,1]\\\cline{2-2}
&X^4+2X^3+11X^2-34X-Y^3-10Y^2+(2X^2-10X-41)Y-32=0\\\cline{2-2}
&25X^{20}-180X^{19}-20X^{18}-1284X^{17}-577X^{16}+18336X^{15}-28656X^{14}\\
&+24644X^{13}-4975X^{12}-71048X^{11}+57344X^9+344064X^8+1089536X^7\\
&+3981312X^6+10518528X^5+36896768X^4+143130624X^3+256901120X^2\\
&+310378496X+402653184\\\cline{2-2}
&^{20}-10X^{19}-9X^{18}+130X^{17}-243X^{16}+5858X^{15}+6995X^{14}-25070X^{13}\\
&+44522X^{12}-14800X^{11}+4096X^9+4096X^8+98304X^7+983040X^6+1843200X^5\\
&+4751360X^4+28934144X^3+54001664X^2+67633152X+113246208\\\cline{2-2}
&-8X^{18}+70X^{17}-100X^{16}+220X^{15}+2148X^{14}-5010X^{13}+2841X^{12}\\
&+952X^{11}+12288X^7+61440X^6+40960X^5+458752X^4+3244032X^3+6291456X^2\\
&+8912896X+12582912\\\cline{2-2}
&X^{14}-2X^{13}-39X^{12}-72X^{11}\\\hline
23&-9+[3,2,5,2]*[5,1,3,1],\quad 4+[3,2,7,2]*[7,1,3,1]\\\cline{2-2}
&X^4+7X^3+27X^2+132X-Y^3+(-8X-6)Y^2+(-9X^2-52X+39)Y+684=0\\\cline{2-2}
&73X^8-546X^7-5395X^6-20918X^5-104397X^4-778894X^3-271080X^2\\
&-1537730X+68337\\\cline{2-2}
&10X^7-572X^6+6588X^5+60922X^4-6432X^3+237609X^2-37183X+146317\\\cline{2-2}
&X^6-506X^5+7017X^4-3420X^3+17834X^2-12608X+2216\\\cline{2-2}
&X-12\\\hline
26&2-[10,2,12,2]*[12,1,10,1],\quad 4+[4,2,9,2]*[9,1,4,1]\\\cline{2-2}
&X^4-7X^3+8X^2-48X-Y^3+(-8X+4)Y^2+(-12X^2-12X-48)Y=0\\\cline{2-2}
&102X^{23}-4274X^{22}+36660X^{21}-88128X^{20}+97312X^{19}-310368X^{18}\\
&-963392X^{17}+293376X^{16}+1023744X^{15}+294656X^{14}-9216X^{13}\\
&+258048X^{10}+5070848X^9+29339648X^8+138493952X^7+511901696X^6\\
&+1364721664X^5+2954887168X^4+4991221760X^3+5486149632X^2\\
&+3288334336X+805306368\\\cline{2-2}
&-12X^{22}-571X^{21}+36436X^{20}-84464X^{19}-234960X^{18}-646160X^{17}\\
&-526784X^{16}+841472X^{15}+1138688X^{14}+303104X^{13}+57344X^9\\
&+3198976X^8+18411520X^7+58720256X^6+182648832X^5+394002432X^4\\
&+401604608X^3+79691776X^2-134217728X-67108864\\\cline{2-2}
&X^{21}-775X^{20}+22904X^{19}-96752X^{18}-104512X^{17}+260144X^{16}\\
&+298752X^{15}+18688X^{14}-35840X^{13}+4096X^8+393216X^7+2600960X^6\\
&+9027584X^5+28573696X^4+67895296X^3+92274688X^2+62914560X+16777216\\\cline{2-2}
&X^{13}(X-12)(X+1)\\\hline
\end{array}
\]
\[
\renewcommand{\arraystretch}{1.2}
\begin{array}{|l|l|l|}\hline
28&-79+[12,2,5,2]*[5,1,12,1],15-[11,2,13,2]*[13,1,11,1]+2*[12,2,5,2]*[5,1,12,1]\\\cline{2-2}
&X^4+288X^3+31097X^2+1491984X+(-5X^2-720X-25920)Y-Y^3+26837568=0\\\cline{2-2}
&102X^{27}+306X^{26}-8442X^{25}-25326X^{24}+206121X^{23}+618363X^{22}\\
&-2018586X^{21}-6055758X^{20}+9733654X^{19}+29200962X^{18}\\
&-27239996X^{17}-81719988X^{16}+48219430X^{15}+144658290X^{14}\\
&-55782076X^{13}-167346228X^{12}+54725846X^{11}+164177538X^{10}\\
&+17711414X^9+53134242X^8+28987609X^7+86962827X^6+926038X^5\\
&+2778114X^4+36358X^3+109074X^2\\\cline{2-2}
&X^{27}+3X^{26}-839X^{25}-2517X^{24}+57267X^{23}+171801X^{22}\\
&-1035053X^{21}-3105159X^{20}+6589699X^{19}+19769097X^{18}\\
&-21809781X^{17}-65429343X^{16}+43462314X^{15}+130386942X^{14}\\
&-55415654X^{13}-166246962X^{12}+46472307X^{11}+139416921X^{10}\\
&-27137621X^9-81412863X^8-605493X^7-1816479X^6-3162117X^5\\
&-9486351X^4+2401X^3+7203X^2-343X-1029\\\cline{2-2}
&-17X^{25}-51X^{24}+2681X^{23}+8043X^{22}-110979X^{21}-332937X^{20}\\
&+1095357X^{19}+3286071X^{18}-4868011X^{17}-14604033X^{16}\\
&+12173210X^{15}+36519630X^{14}-18835474X^{13}-56506422X^{12}\\
&+18755674X^{11}+56267022X^{10}-9270443X^9-27811329X^8+11672045X^7\\
&+35016135X^6+1771693X^5+5315079X^4+197225X^3+591675X^2-49X-147\\\cline{2-2}
&X^2(X^2-49)(X^2-1)^7(X+3)\\\hline
29&3-[13,2,11,2]*[11,1,13,1],-4+[4,2,11,2]*[11,1,4,1]-2[13,2,11,2]*[11,1,13,1]\\\cline{2-2}
&X^4+7X^3+16X^2-10X+(-4X+2)Y^2-Y^3+144+(4X^2-10X+35)Y=0\\\cline{2-2}
&31X^{10}+9X^9-159X^8+205X^7+9501X^6-50799X^5+165072X^4-324935X^3\\
&+440497X^2-350079X+137927\\\cline{2-2}
&7X^9+104X^8-1134X^7+6377X^6-25276X^5+66921X^4-127602X^3+161922X^2\\
&-131261X+49677\\\cline{2-2}
&X^8-140X^7+1128X^6-4832X^5+13228X^4-24294X^3+29790X^2-22321X+7948\\\cline{2-2}
&X-6\\\hline
31&10+4,2,9,2]*[9,1,4,1],-22+[3,2,5,2]*[5,1,3,1]+[4,2,9,2]*[9,1,4,1]\\\cline{2-2}
&X^4-21X^3+235X^2-3061X+10304-Y^3-(4X+24)Y^2\\
&+(4X^2-210X+409)Y=0\\\cline{2-2}
&-18X^{11}+4997X^{10}-196897X^9+2788515X^8-18906271X^7+72162582X^6\\
&-248096148X^5+1325621567X^4-6136101579X^3+16559053779X^2\\
&-23236346958X+13294184235\\\cline{2-2}
&X^{10}-461X^9-26628X^8+1849583X^7-31107208X^6+241112657X^5\\
&-995412091X^4+2185571830X^3-2080251619X^2-266712928X+1365873417\\\cline{2-2}
&217X^8-37014X^7+1334457X^6-19397444X^5+145182734X^4-615528808X^3\\
&+1500491171X^2-1968040952X+1079470840\\\cline{2-2}
&X-8\\\hline
\end{array}
\]
}
For other cases we list only functions $F_1,F_2,\dots,F_{g_0(N)+1}$.\vspace{3mm}\newline
(iv) the case $g_0(N)=2$ and $37\leq N\leq 50$.
{\scriptsize
\[
\renewcommand{\arraystretch}{1.2}
\begin{array}{|l|l|l|l|}\hline
N&F_1&F_2&F_3\\\hline
37&-\frac 12[6,4,13,4]*[13,1,6,1]&-[13,2,17,2]*[17,1,13,1]&[3,2,5,2]*[5,1,3,1]\\\hline50&-\frac 83+\frac 13(-h_1+h_2-h_3+2h_4)&-2+h_1+h_3-h_4&3-h_1\\\hline
\end{array}
\]
}
Here $h_1=[10,5,15,5]*[15,1,10,1],h_2=[10,8,24,8]*[24,1,10,1]$,$h_3=[22,7,3,7]*[3,1,22,1]$,$h_4=[10,4,24,4]*[5,1,16,1]$.\vspace{3mm}\newline
(v) The case $g_0(N)=3$.
{\scriptsize
\[
\renewcommand{\arraystretch}{1.2}
\begin{array}{|l|l|l|}\hline
N&F_1 & F_2\\\cline{2-3}
 &F_3&F_4\\\hline
30&-[4,2,14,2]*[14,1,4,1] &-[10,5,9,5]*[9,1,10,1]\\\cline{2-3}
 &-[4,2,11,2]*[11,1,4,1] &-[13,2,7,2]*[7,1,13,1]\\\hline
33& [16,2,5,2]*[5,1,16,1]&[13,2,11,2]*[11,1,13,1]\\\cline{2-3}
&[4,2,15,2]*[15,1,4,1]&[4,2,7,2]*[7,1,4,1]\\\hline
34&-[3,5,14,5]*[14,1,3,1]&[4,2,16,2]*[16,1,4,1]\\\cline{2-3}
&-[4,2,13,2]*[13,1,4,1]&[6,2,11,2]*[11,1,6,1]\\\hline
35&-[9,2,16,2]*[16,1,9,1]-[7,2,15,2]*[15,1,7,1]&[9,2,16,2]*[16,1,9,1]\\\cline{2-3}
&[9,2,14,2]*[14,1,9,1]+[9,2,16,2]*[16,1,9,1]&[3,2,17,2]*[17,1,3,1]\\
&+[7,2,15,2]*[15,1,7,1]& \\\hline
39&[18,2,14,2]*[14,1,18,1]&-[15,2,3,2]*[3,1,15,1]\\
&&-F_4-3F_3-3F_1\\\cline{2-3}
&[4,2,9,2]*[9,1,4,1]-[18,2,14,2]*[14,1,18,1]&[16,2,14,2]*[14,1,16,1]+2F_3+F_1\\\hline
40&\frac 12[17,3,7,3]*[7,1,17,1]&[6,4,14,4]*[14,1,6,1]\\\cline{2-3}
&-\frac 12[5,3,19,3]*[19,1,5,1]&-[2,3,18,3]*[18,1,2,1]\\\hline
41&-[7,2,5,2]*[5,1,7,1]&\frac12[3,2,5,2]*[5,1,3,1]\\\cline{2-3}
&-[4,2,9,2]*[9,1,4,1]-F_1-F_2&-[4,2,12,2]*[12,1,4,1]-F_3+2F_2\\\hline
43&-\frac 12[3,4,5,4]*[5,1,3,1]&[12,2,5,2]*[5,1,12,1]\\\cline{2-3}
&[7,2,5,2]*[5,1,7,1]&[4,2,9,2]*[9,1,4,1]\\\hline
45&\frac 13[7,8,15,8]*[15,1,7,1]&-\frac 12[6,4,21,4]*[21,1,6,1]\\\cline{2-3}
&\frac 13[12,3,21,3]&[6,3,21,3]*[21,1,6,1]\\\hline
48&\frac 14[20,4,8,4]&[11,17,8,17]*[8,1,11,1]\\\cline{2-3}
&\frac 12[13,5,3,5]*[3,1,13,1]&-[13,2,6,2]*[6,1,13,1]\\\hline
 \end{array}
 \]
}
v) The case $g_0(N)=4$.
{\scriptsize
\[
\renewcommand{\arraystretch}{1.2}
\begin{array}{|l|l|}\hline
N&F_1,~ F_2,~F_3,~F_4,~F_5\\\hline
38&-[4,2,18,2]*[18,1,4,1]+[18,1,2,1]*[5,1,5,2]+[6,2,13,2]*[13,1,6,1],~[18,1,2,1],\\
&[4,2,18,2]*[18,1,4,1],~-[14,4,5,4],~-[15,2,4,2]\\\hline
44&-[11,4,10,4]*[10,1,11,1],~\frac 12[20,4,2,4],~[21,1,2,1],~[17,20,5,20]*[5,1,17,1],\\&-\frac 12[2,11,20,11]\\\hline
47&-[9,23,2,23]*[2,1,9,1]+[13,20,2,20]*[2,1,13,1]+[1,2,5,2],~-[1,2,5,2],\\
&-[9,23,2,23]*[2,1,9,1],~[5,1,2,1],~[5,1,3,1]\\\hline
\end{array}
\]
}
(vi) The case $g_0(N)=5$.
{\scriptsize
\[
\renewcommand{\arraystretch}{1.2}
\begin{array}{|l|l|}\hline
N&F_1,~ F_2,~F_3,~F_4,~F_5,~F_6\\\hline
42&[19,2,5,2]*[5,1,19,1],~\frac 12[7,6,14,6],~[17,11,4,11]*[4,1,17,1],\\
&[16,2,5,2]*[5,1,16,1],~[18,3,13,3]*[13,1,18,1],~[10,4,11,4]*[11,1,10,1]\\\hline
46&\frac 12+[19,8,3,8]+\frac 12[16,14,1,14]+[1,18,2,18]*[2,1,20,1]+[8,3,2,3]*[2,1,8,1],\\&-\frac 32[1,7,2,7]*[2,1,5,1],-[22,1,2,1],~-[4,2,22,2]*[22,1,4,1],~[20,2,3,2]*[3,1,20,1],\\
&-[15,1,8,1],~-[19,2,4,2]\\\hline
\end{array}
\]
}
{\bf Acknowledgement.} The author wishes to thank Naoya Nakazawa for his assistance in computing some cases in our tables. 
\vspace{5mm}
{\small
Faculty of Liberal arts and Sciences\newline
Osaka Prefecture University \newline 
1-1 Gakuen-cho, Sakai, Osaka\\
 599-8531 Japan\newline
e-mail:\quad ishii@las.osakafu-u.ac.jp }
\end{document}